\documentclass[11pt]{article}
\pdfoutput=1

\usepackage[T1]{fontenc}
\usepackage[utf8]{inputenc}
\usepackage{lmodern}
\usepackage{microtype}

\usepackage{amsmath,amssymb,amsthm,latexsym,mathtools}
\usepackage{graphicx}
\usepackage{xcolor}
\usepackage{booktabs}
\usepackage{nicematrix}
\usepackage[ruled,vlined]{algorithm2e}
\usepackage{enumitem}
\usepackage[left=1in,right=1in,top=0.8in,bottom=0.9in]{geometry}
\usepackage{changepage}
\usepackage{hyperref}
\usepackage{cleveref}

\usepackage[numbers,sort&compress]{natbib}

\newcommand{\hl}[1]{#1}
\newcommand{\highlighting}[1]{#1}

\newcommand{\Ex}{\ensuremath{\mathbb{E}}}

\def\maxim{\mathop{\textup{maximize}}}

\providecommand{\natexlab}[1]{#1}

\theoremstyle{plain}

\theoremstyle{definition}

\title{Markovian restless bandits and index policies: A review}
\author{Jos\'e Ni\~no-Mora\\
Departamento de Estad\'{\i}stica, Universidad Carlos III de Madrid\\
28903 Getafe (Madrid), Spain\\
\texttt{jose.nino@uc3m.es}}
\date{} 

\begin{document}

\maketitle

\begin{abstract}
The restless multi-armed bandit problem is a paradigmatic modeling framework for optimal dynamic priority allocation in  stochastic models of wide-ranging applications that has been widely investigated and applied since its inception in a seminal paper by Whittle in the late 1980s. 
The problem has generated a vast and fast-growing literature from which a significant sample is thematically organized
and reviewed in this paper.
While the main focus is on priority-index policies due to their intuitive appeal, tractability, asymptotic optimality properties, and often strong empirical performance, other lines of work are also reviewed.  Theoretical and algorithmic developments are discussed, along with diverse applications. The main goals are to highlight the remarkable breadth of work that has been carried out on the topic and to stimulate further research in the field.
\end{abstract}



\medskip
\noindent\textbf{Note:} Published in \emph{Mathematics} \textbf{11}, 1639 (2023). DOI:\ \url{https://doi.org/10.3390/math11071639}

\section{Introduction}
\label{s:intro}
In the much-studied \emph{\hl{multi-armed bandit problem} 
} (MABP), one is presented with a number of generic dynamic and stochastic \emph{\hl{projects}} modeled as binary-action (active/selected or passive/rested) \emph{Markov decision processes} (MDPs) (see, e.g.,~Puterman~\cite{put94}) with the aim of maximizing the expected reward accrued by dynamically selecting one project to engage at each time. 
The~problem, named after 
a
multi-armed slot machine (or ``bandit'') ---whence the projects are commonly referred to in the literature as \emph{\hl{bandits}} or \emph{\hl{arms}}--- is widely regarded as a major modeling framework for addressing the \emph{exploration vs.\ exploitation trade-off} in widely diverse settings. A key assumption is that project states remain frozen while passive, which renders the problem tractable in some relevant scenarios, most notably under the expected total infinite-horizon (geometrically) discounted-reward criterion. In this case, as first shown by Gittins and Jones~\cite{gijo74}, one can attach to each project a scalar function of its state called an \emph{index} that is based solely on the project's characteristics, such that the resulting \emph{priority-index policy}, which selects at each time a project with highest index, is optimal. 

In a seminal work, Whittle~\cite{whit88b} first considered a model that dropped the assumption of frozen states by allowing rested projects to change states. Projects then become \emph{restless}, giving rise to the \emph{restless multi-armed bandit problem} (RMABP). Whittle outlined potential applications to demonstrate the expanded modeling power of the new modeling framework, and, to overcome its intractability due to the \emph{curse of dimensionality}, he proposed a tractable heuristic priority-index policy based on Lagrangian relaxation and decomposition ideas. He further conjectured a form of asymptotic optimality for this \emph{Whittle index policy}. 

The impact of Whittle's work took off to a slow start, but~over the following decade researchers began to realize the potential of Whittle's model and solution approach. The~number of researchers interested in the RMABP and its variants, as well as~the resulting number of published papers on the subject, has grown steadily since then and has picked up steam in recent years. The literature on the RMABP, whether on its theoretical, algorithmic, or~application aspects, is currently vast to~the point where it is virtually infeasible for researchers to keep up to date with the latest advances in the field. 

This is the motivation of the present paper, which reviews and highlights a hopefully representative albeit limited and incomplete sample of the enormous body of work published on the topic by a myriad of researchers over the thirty-plus years since its inception by Whittle, emphasizing breadth over depth. The main goal is to stimulate further research on this subject both by bringing it to the attention of researchers who may not have  encountered it previously, and by providing an overview of its wide-ranging possibilities and open problems for those who may have missed important research directions. This review focuses on solution approaches that develop, test, and analyze index policies due to their intuitive appeal, tractability, asymptotic optimality properties, and often strong empirical performance. The author has strived to go well beyond his comfort zone to provide a balanced coverage of the subject, minimizing bias due to his familiarity with index policies for models with known parameters, although some inevitable bias may remain.

The review is organized as follows. 
Section~\ref{s:antec} surveys the antecedents to the RMABP, in~particular, the classic MABP and the Gittins index policy. 
Section~\ref{s:wrmabpip} formulates the RMABP and outlines the Whittle index policy.
Section~\ref{s:cca} reviews works on the complexity, approximation, and~relaxations of the RMABP.
Section~\ref{s:api} focuses on indexability, that is, the existence of the Whittle index and extensions.
Section~\ref{s:ic} discusses works on means of computing the Whittle index.
Section~\ref{s:oip} considers works that establish the optimality of the myopic index policy for the RMABP, whereas 
Section~\ref{s:aowip} reviews the asymptotic optimality of index policies.
Section~\ref{s:marb} surveys multi-action restless bandits.
Section~\ref{s:lagpol} considers policies that are different from Whittle's and based on Lagrangian and fluid relaxations. 
Section~\ref{s:reinflearn} addresses reinforcement learning and, in~particular, Q-learning solution approaches.
Section~\ref{s:learning} surveys works on the RMABP from the perspective of online learning.
Sections~\ref{s:appl} and \ref{s:applPO} are devoted to works on applications of the RMABP in diverse settings. The~former section focuses on MDP models and the latter on \emph{partially observable MDP} (POMDP) models.
Finally, Section~\ref{s:concl} concludes this paper.

\section{Antecedents: Multi-Armed Bandits and the Gittins Index~Policy}
\label{s:antec}
The investigation of bandit-like problems has its early roots in the seminal work of Thompson~\cite{thompson33}, who outlined a research-planning problem that aimed to allocate (``apportion'' in the words of Thompson) two Bernoulli treatments with unknown success probabilities to a set of individuals based on prior statistical evidence using Bayes' rule to~reduce the expected number that received inferior treatment. Thompson \cite{thompson35} further clarified the use and implementation of his proposed approach for dynamic sequential allocation (the basis of the widely investigated and applied \emph{Thompson sampling} method), sketched an extension to multiple Bernoulli treatments, and~even reported the results of a pioneering simulation study avant la lettre, illustrating its practical~effectiveness.

This precursor to the MABP was later refined and cast into the framework of the then-budding theory of the \emph{sequential design of experiments} (see Wald~\cite{wald47} for its first systematic account), which addressed the design of sequential sampling procedures ``in which the size and composition of the samples are not fixed in advance but are functions of the observations themselves'', as~stated by Robbins~\cite{robbins52}. 
In that paper, Robbins considered the two-armed bandit problem of sampling from two Bernoulli populations to maximize the expected number of successes from a fixed number of observations under the practically relevant assumption that population parameters are unknown. He thus adopted a \emph{learning} perspective and showed that one can design sampling policies that, as the number of observations grows, asymptotically attain the maximum expected reward per observation.

Bradt~et~al.~\cite{bjk56} adopted a Bayesian perspective and considered a finite-horizon two-armed bandit problem to maximize the expected sum of a fixed number of observations through sequential sampling from two unknown statistical populations. Their approach assumed a known prior distribution. 
They identified conditions under which the myopic index policy is optimal but, in general, noted the difficulty of elucidating the structure of optimal policies. Yet, in the case of two Bernoulli populations, one with known parameter, they showed that the optimal policy is characterized by  a critical quantity, which we would call an \emph{index}, a function of both the remaining number of observations and the current posterior distribution of the unknown population. At each time, it is optimal to sample from the unknown population if and only if its current index value is greater than or equal to the known population's parameter. Bellman~\cite{bellman56} considered the problem under the infinite-horizon discounted reward criterion, also establishing the optimality of an \emph{index policy}. In this case, the index is a function only of the posterior distribution of the success probability of the unknown population.

MABPs were later cast into the more general framework of multi-stage decision processes, in~particular, MDPs (see the pioneering monograph of Bellman (\cite{bellman57} Chapter 11) and the later work of Howard~\cite{howard60}, and for~more recent accounts of MDPs, see, e.g.,~ the comprehensive textbooks  of Puterman~\cite{put94} and Bertsekas~\cite{bertsek17,bertsek12}).
 
 A general formulation of a standard version of the MABP in the MDP setting is as follows. A decision maker is faced with a finite collection of $N$ reward-yielding dynamic and stochastic \emph{projects}. The \emph{state} $X_{n, t}$ of project $n = 1, \ldots, N$ evolves over discrete time periods $t = 0, 1, 2, \ldots$ over an infinite horizon across the \emph{state space} $\mathcal{X}_n$. The evolution and reward of project $n$ in period $t$ depend both  on the current state $X_{n, t}$ and the chosen \emph{action} $A_{n, t}$, which can take two values: $1$ (\emph{active}, i.e.,~\emph{engaging} the project) and $0$ (\emph{passive}, i.e.,~\emph{resting} it). When the project occupies state $X_{n, t}= i_n$ and action $A_{n, t} = a_n$ is selected, the~project yields an expected reward $r_n^{a_n}(i_n)$, which is time-discounted with the factor $0 < \beta < 1$ and~its state moves to $X_{n, t+1}= j_n$ with a transition probability of $p_n^{a_n}(i_n, j_n)$, where $p_n^{0}(i_n, j_n) = \delta_{i_n j_n}$ (Kronecker's delta) so passive projects do not change state. It is usually assumed that $r_n^{0}(i_n) \equiv 0$ for all $i_n \in \mathcal{X}_n$, i.e.,~passive projects give no rewards but~this assumption is nonessential. At each time $t$, the~decision maker observes the \emph{joint state} \highlighting{$\mathbf{X}_t = (X_{n, t})_{n=1}^N$} and then selects one project so the \emph{joint action} \highlighting{$\mathbf{A}_t = (A_{n, t})_{n=1}^N \in \{0, 1\}^N$} must~satisfy 
\begin{equation}
\label{eq:mabpAlteqs1}
\sum_{n = 1}^N A_{n, t} = 1, \quad t = 0, 1, 2, \ldots
\end{equation}

The class of non-anticipative \emph{policies} that prescribes these feasible action choices is denoted as $\Pi$  and the expectation starting from the joint state $\mathbf{i}$ under policy $\pi$ is denoted as $\Ex_{\mathbf{i}}^\pi[\cdot]$. The goal of the \emph{infinite-horizon discounted MABP} is to find an \emph{optimal policy} $\pi^* \in \Pi$, which, for~\emph{any} initial joint state $\mathbf{i}$, maximizes the expected total discounted reward earned over an infinite horizon, which is given by
\begin{equation}
\label{eq:mabpobj}
\Ex_{\mathbf{i}}^\pi\bigg[\sum_{t = 0}^\infty \sum_{n = 1}^N r_n^{A_{n, t}}(X_{n, t}) \beta^t\bigg],
\end{equation}

Note that applying the conventional solution approach based on solving numerically Bellman's optimality equations using \emph{dynamic programming} (DP) is hindered by the \emph{curse of dimensionality} as, even in the finite state case, the number of joint states and hence of equations increases exponentially with the number of~projects.

After being considered intractable for a long time, Gittins and Jones~\cite{gijo74} made a breakthrough by demonstrating that the optimal policies for the above MABP have a strikingly simple structure. Specifically, for each project $n$, there exists an \emph{index} $\varphi_n(i_n)$, which is a function of the project state $i_n$, that is based solely on its characteristics (rewards and transition law) such that the resulting \emph{index policy}, which selects a project with the highest index value at each time, is optimal. This index was introduced in the aforementioned work of Bellman~\cite{bellman56} for a special Bernoulli bandit model (see also the paper by Gittins~\cite{gi79} and the monographs by Gittins~\cite{gi89} and Gittins ~et~al.~\cite{gigwe11}).

Researchers have utilized various approaches to offer different proofs for these celebrated results on the optimality of what is widely known as the \emph{Gittins index} policy (see, e.g.,~Whittle~\cite{whittle80}, Varaiya~et~al.~\cite{vawabu85}, Weber~\cite{weber92}, Tsitsiklis~\cite{tsitsik94}, and Bertsimas and Ni\~{n}o{-}Mora~\cite{bnm96}). Furthermore, the results have been extended to incorporate variations that preserve the optimality of index policies, e.g.,~random project arrivals (Whittle~\cite{whit81}), non-Markovian evolution (Varaiya~et~al.~\cite{vawabu85}), the minimization of the expected cost for one project to reach a target state (Dumitriu~et~al.~\cite{dumitriuetal03}), and constrained arm switches (Bao~et~al.~\cite{baoetal21}), among others. 
 
In independent ground-breaking work, Klimov~\cite{kl74} presented a landmark result analogous to that in~\cite{gijo74} by first establishing the optimality of an index policy for scheduling a multi-class M/G/1 queue with Bernoulli job feedback and linear holding costs under the average criterion, which is a semi-Markov MABP with project arrivals. The results were extended to the more general \emph{branching bandit} problem by
Meilijson and Weiss~\cite{meilijWeiss77}, who analyzed it under the average criterion, and Weiss~\cite{weiss88},  who analyzed it under the discounted criterion. For~more on Klimov's results, the~reader is referred to the review article by Ni\~{n}o{-}Mora~\cite{nmKlimovsEORMS}. 

\section{Restless Multi-Armed Bandits and the Whittle Index~Policy}
\label{s:wrmabpip}
Whittle~\cite{whit88b} introduced a major extension to the classic MABP by allowing  rested projects to change state, in which case they are called \emph{restless bandits}, with state transitions being independent across projects. The model was also considered in independent work by O'Flaherty~\cite{oflaherty89}.
Whittle noted that the restless feature substantially expanded the modeling power, which, however, came at the expense of tractability, as~index policies are generally not optimal for the resulting RMABP.
Whittle further allowed the number of projects selected at each time, $M$,  to not necessarily be one. 
The \emph{joint action} $\mathbf{A}_t = (A_{n, t})_{n=1}^N \in \{0, 1\}^N$ must thus satisfy
\begin{equation}
\label{eq:mabpAlteqsK}
\sum_{n = 1}^N A_{n, t} = M, \quad t = 0, 1, 2, \ldots
\end{equation}

Whittle focused on the average-reward criterion, the goal being to find a policy satisfying (\ref{eq:mabpAlteqsK}) that maximizes the steady-state reward
\begin{equation}
\label{eq:mabpobjta}
\Ex^\pi\bigg[\sum_{n = 1}^N r_n^{A_n}(X_n)\bigg],
\end{equation}
where $(X_n, A_n)$ has the steady-state distribution of state-action pairs for project $n$ under policy $\pi$ and~$\Ex^\pi[\cdot]$ denotes the expectation under this policy, which does not depend on the initial joint state under suitable regularity conditions. He considered a \emph{relaxation} of this problem, where the sample-path constraint (\ref{eq:mabpAlteqsK}) is \emph{relaxed} by requiring only that $M$ projects be active \emph{on average} so that
\begin{equation}
\label{eq:Kpaoa}
\Ex^\pi\bigg[\sum_{n = 1}^N A_n\bigg] = M.
\end{equation}

Whittle showed that this relaxed problem, whose optimal value provides an upper bound on that of the original problem, was amenable to a Lagrangian solution approach by attaching a Lagrange multiplier $\lambda$  to constraint (\ref{eq:Kpaoa}). In the resulting \emph{Lagrangian relaxation}, the constraint is brought into the objective, with the multiplier playing the economic role of a \emph{subsidy for passivity}. He further noted that the Lagrangian relaxation allows a decomposition, which entails single-project subproblems of the form
\begin{equation}
\label{eq:spspwhit}
\maxim_{\pi_n \in \Pi_n} \, \Ex^\pi\big[r_n^{A_n}(X_n) + \lambda (1-A_n)\big],
\end{equation}
where $\Pi_n$ is the class of admissible policies for operating project $n$ \emph{in isolation}. Now, on the one hand, one can consider the \emph{Lagrangian dual problem}, which is to find the best multiplier $\lambda$ for which the Lagrangian relaxation yields the best possible upper bound. On the other hand, one must consider the properties of optimal policies for the subproblems in (\ref{eq:spspwhit}).

Whittle proposed the following \emph{indexability} property for a single project: as the passivity subsidy  $\lambda$ increases from $-\infty$ to $+\infty$, the~set of project states where it is optimal in (\ref{eq:spspwhit}) to rest the project increases monotonically from the empty set  to the full project state space. The index $\varphi_n(i_n)$ for project $n$ in state $i_n$ is then defined as the value of $\lambda$ at which both actions are optimal in $i_n$, which is unique under~indexability. 

If each constituent project in an RMABP is indexable, each single-project subproblem in (\ref{eq:spspwhit}) is solved optimally by an index policy with index $\varphi_n(i_n)$. It is optimal to select (resp.\ rest) project $n$ in state $i_n$ if and only if $\varphi_n(i_n) \geqslant \lambda$ (resp.\ $\varphi_n(i_n) \leqslant \lambda$). In such a case, Whittle proposed to use the resulting priority-index policy as a heuristic rule for the RMABP. At each time, select $M$ projects with the highest index values. He further conjectured that, although the \emph{Whittle index policy} will generally be suboptimal, it should have a form of \emph{asymptotic optimality}. In the case of a population of projects with a fixed number of distinct project types in fixed proportions, as $M$ and $N$ grow to infinity with a fixed ratio $M/N$, the average reward per project should converge to that under an optimal policy.

Whittle also proposed to use this approach to define an index policy for the discounted RMABP, noting that for non-restless projects his index reduces to the Gittins index.

\section{Complexity, Approximation, and~Relaxations}
\label{s:cca}
In contrast to the classic MABP for which the Gittins index policy allows breaking the curse of dimensionality under the infinite-horizon discounted-reward criterion, the RMAB is, in general, computationally intractable, as first shown by Papadimitriou  and Tsitsiklis~\cite{papTsik99}. Among a host of ground-breaking complexity results, they showed the RMABP to be PSPACE-hard, even for projects with deterministic state transitions. As they stated, this  ``is considered a much more convincing evidence of intractability than~NP-hardness.''

Regarding theoretical performance guarantees for heuristic policies, Guha~et~al.~\cite{guhaetal10} presented the first \emph{constant-factor approximation algorithms} based on \emph{linear programming} (LP) duality for several classes of RMABPs under the average-reward criterion, including a $(2~+~\epsilon$)-approximation algorithm for the  \emph{feedback RMABP}, which includes the \emph{partially observable MDP} (POMDP) RMABP model for opportunistic spectrum access considered in Liu and Zhao~\cite{liuZhao08,liuZhao10} and in Le Ny~et~al.\  \cite{lenyetal08}. They considered an index policy that is closely related to Whittle's, and extended their results to the wider class of \emph{monotone bandits}.

Other researchers have followed in the footsteps of~\cite{guhaetal10} to design constant-factor approximation algorithms for related RMABP classes that incorporate additional relevant features (see, e.g.,~ Wan and Xu~\cite{wanXu15}, Xu and Song~\cite{xuSong16}, and~Xu and Wang~\cite{xuWang20}, which incorporate, respectively,  project weights, channel interferences, and~multiple~constraints).

As outlined above, Whittle proposed a problem relaxation to approximate the RMABP. This idea was extended by Bertsimas and Ni\~{n}o{-}Mora~\cite{benm00} who considered a hierarchy of $N$ increasingly stronger \emph{linear programming} (LP) relaxations (where $N$ is the number of projects), with increasing size and hence complexity. The first relaxation in the hierarchy is that in~\cite{whit88b}, whereas the last is the standard exact formulation of~exponential size  in $N$. They further proposed and tested an alternative index policy based on  the optimal solution to the first-order relaxation, which does not require projects to be~indexable.

\textls[-10]{Whittle's relaxation approach to the RMABP was further extended by
Hawkins~\cite{hawkins03}} and then by Adelman  and Mersereau~\cite{adelMerse08} into the broader setting of so-called \emph{weakly coupled} MDPs, in which a set of sample-path linking constraints couples a collection of otherwise independent constituent subproblems. They compared bounds and policies obtained from two different problem relaxations, which were based, respectively, on LP-based \emph{approximate DP} (see, e.g.,~the monographs by Powell~\cite{powell11} and Bertsekas~\cite{bertsek12}) and Lagrangian relaxation. Among their results, they  showed that both relaxations entail fitting an additively separable value function approximation to the Bellman optimality equations and established that the approximate DP bound was tighter than the one based on Lagrangian~relaxation.

The analyses of relaxations in~\cite{hawkins03,adelMerse08} were further refined by
Brown and Zhang~\cite{brownZhang22}, who  provided theoretical justification in~the form of sufficient conditions for~the empirically observed fact that the gap between the two different upper bounds on the optimal value considered in~\cite{adelMerse08} was typically~small.

\section{Indexability}
\label{s:api}
Whittle~\cite{whit88b} pointed out that ``One would very much like to have simple sufficient conditions for indexability; at the moment, none are known.'' Most works where Whittle's index policy has been applied to particular RMABP models (see, e.g.,~ the references in Sections~\ref{s:appl} and \ref{s:applPO})
have followed the following scheme to argue that individual projects are indexable: (i) show that the subproblems in (\ref{eq:spspwhit}) (or their discounted counterparts) can be solved optimally using a certain family of structured policies for any value of the multiplier $\lambda$, typically threshold policies for linearly-ordered project state spaces; (ii) obtain an optimal policy within such a family of policies as a function of $\lambda$ for these subproblems, e.g.,~an optimal threshold policy; and (iii) show that the function derived in the previous step satisfies a required form of monotonicity in $\lambda$, under which the Whittle index is then obtained by inverting this function. In some works, step (i) is carried out rigorously, e.g., using threshold policies in Le Ny~et~al.~\cite{lenyetal08}, Liu and Zhao~\cite{liuZhao08,liuZhao10}, and~Liu~et~al.~\cite{liuWebZhao11}; using stopping rules in~Fryer and Harms~\cite{fryerHarms18}; and using another family of policies in Caro and Yoo~\cite{caroYoo10}. Yet, often, the required optimality of the family of structured policies under consideration is merely postulated or conjectured as researchers either bypass it, as, e.g.,~ in  Whittle (\cite{whit96} Chapter 14.6) and in Veatch and Wein~\cite{vewe96}, or find that the proof is elusive, as in~Dance and Silander~\cite{danceSi15}.

A different approach based on \emph{partial conservation laws} (PCLs) that provides sufficient indexability conditions for general restless bandits has been developed by the author in a series of papers \cite{nmaap01,nmmp02,nmmor06,nmtop07,nmmor20}. It proceeds by (i) positing a family of structured policies that one guesses might be optimal for the single-project subproblems in (\ref{eq:spspwhit}); (ii) proving that a certain \emph{marginal work} project metric is positive for the postulated family of policies; and (iii) showing that the \emph{marginal productivity index} computed by an adaptive-greedy algorithm given in~\cite{nmaap01} is generated in a monotonic fashion. Steps (ii) and (iii) can be carried out either numerically or analytically. The main result is that under the conditions in (ii) and (iii), it follows simultaneously in one fell swoop that both the postulated family of policies is optimal for the single-project subproblems \emph{and} the obtained marginal productivity index is the project's Whittle index. Furthermore, under the condition in (ii), the model is indexable consistently with the postulated family of policies if and only if the condition in (iii) holds.
This approach allowed to overcome the long-standing problem of establishing the optimality of threshold policies and proving indexability for the scalar Kalman filter restless bandit model, as demonstrated in the groundbreaking work of Dance and Silander~\cite{danceSi19}.

\section{Whittle Index~Computation}
\label{s:ic}
\textls[-10]{In some restless bandit models, the Whittle index can be derived in closed form, e.g., in Whittle~\cite{whit88b}, Whittle (\cite{whit96} Chapter 14.6), Veatch and Wein~\cite{vewe96}, and Liu~et~al.~\cite{liuWebZhao11}. For other models where the Whittle index needs to be computed, an~exact algorithm applied to a general $n$-state restless bandit satisfying PCLs was provided in~\cite{nmaap01,nmmp02}} and~a fast-pivoting implementation with complexity $O(n^3)$ was presented in~\cite{jnmsmct07b,nmmath20}, extending the author's work in~\cite{nmijoc07} on efficient computation of the Gittins index. Ref.~\cite{jnmsmct07b} further provided an index algorithm, also with $O(n^3)$ complexity but with a larger leading constant, for checking indexability and computing the index without the need to satisfy the PCLs. Note that these algorithms \emph{both} check for indexability \emph{and} compute the Whittle index when it exists. In special cases, the complexity is reduced, e.g.,~in the birth--death queueing model in~\cite{nmmp02}, where it is shown that threshold policies are optimal and the Whittle index is computed in $O(n)$ time. For models where projects have a continuous, real state space, the~PCL approach also provides a means of computing the Whittle index, as shown in~\cite{nmmor20}.

Another approach to computation of the Whittle index  was provided by Qian~et~al.~\cite{qianetal16}, who   introduced sufficient conditions for indexability in their model of optimal patrol policies and provided an algorithm to test for indexability. Given indexability, they proposed to use binary search to compute the Whittle~index.
 
 Akbarzadehand and Mahajan~\cite{akbarMahaj22} extended the adaptive-greedy algorithm in~\cite{nmaap01,nmmp02} to a version that can compute the Whittle index for an arbitrary indexable restless bandit. They further presented an efficient implementation of their algorithm with $O(n^3)$ complexity for an $n$-state project and provided alternative sufficient conditions for indexability. 

\section{Optimality of the Myopic~Policy}
\label{s:oip}
Although index policies are generally suboptimal for the RMABP, a substantial amount of research has been devoted to identifying conditions under which the \emph{myopic policy}, which prioritizes projects based on their current expected reward when selected, is optimal. However, these results typically require strong assumptions about the constituent projects, most notably that they are homogenous.

The first work in this area was by O'Flaherty~\cite{oflaherty89}, who provided sufficient conditions for the myopic policy to be optimal in two-project RMABPs. The paper outlined ideas supporting the validity of these conditions but did not provide formal proofs.

Ehsan and Liu~\cite{ehsanLiu09} considered an RMABP model for optimal dynamic server allocation in a multi-class single-server discrete-time queue with delayed backlog information and convex nondecreasing holding costs to minimize  the expected total discounted cost over a finite or infinite horizon. They derived sufficient conditions for a myopic policy to be optimal in certain regions of the state space, in the cases of batch and Poisson job arrivals, which take  the form of sufficient separation among indices. 

Ahmad~et~al.~\cite{ahmadetal09} considered a POMDP RMABP model with  $M = 1$ active projects at each time for dynamic multi-channel access with identical channels and~demonstrated the optimality of the myopic policy under both discounted- and average-reward criteria, provided that state transitions exhibit a positive correlation over time. Liu~et~al.~\cite{liuZhaoKrishna10} incorporated  error-prone channel state detection and demonstrated the optimality of the myopic policy for the two-channel~case.

In a series of papers, Wang and coworkers investigated conditions under which the myopic policy is optimal for several POMDP RMABP models of optimal dynamic multi-channel access. In \cite{wangChen12a}, Wang~et~al. considered the dynamic multi-channel access model with identical channels and the class of so-called standard reward functions, for~which  a closed-form condition on the discount factor ensures optimality of the myopic policy under the discounted criterion. The results were also extended to the average-reward criterion. Wang and Chen~\cite{wangChen12} considered a general setting and introduced three axioms that characterized the so-called regular functions, which yielded closed-form conditions under which the myopic policy is optimal. The results were demonstrated in RMAB problems arising in multi-channel opportunistic access with heterogeneous channels. These results were extended by 
Wang~et~al.~\cite{wangetal14} to a framework based on the so-called g-regular functions.
Wang~et~al.~\cite{wangetal16} considered a multi-channel access problem where each of $N$ channels is modeled as a multi-state Markov chain rather than as a two-state Markov chain as in prior work. They demonstrated that the myopic policy reduces channel selection to a round-robin policy whose optimality is established for accessing $M = N - 1$ out of $N$ channels. Wang~et~al.~\cite{wangetal22} further extended this work by  identifying two sets of sufficient conditions on both the eigenvalues and eigenmatrices resulting from channel-state transition matrices that guarantee the optimality of the myopic policy (see also the monograph by Wang and Chen~\cite{wangChen21}).

Ouyang and Teneketzis~\cite{ouyangTenek14} considered a POMDP dynamic multi-channel access RMABP model, where the underlying channel-state Markov chain has an arbitrary finite number of states.  They presented sufficient conditions  for the optimality of a myopic sensing policy over a finite-time horizon under discounted and undiscounted-reward criteria.

Blasco and G\"und\"uz~\cite{blascoGun15} proposed a POMDP RMABP model for a multi-access wireless network with transmitting nodes, each with an energy-harvesting (EH) device and a finite-capacity rechargeable battery, with the goal of maximizing throughput. Under certain conditions on the EH processes and battery sizes, the optimality of the myopic policy is shown. 

The Age of Information (AoI) is a popular metric for capturing information freshness, based on the time elapsed since the most recently delivered packet in a communication node. Kadota~et~al.~\cite{kadotaetal18} proposed a model for minimizing the AoI and demonstrated that, in~symmetric networks, the~myopic policy, which prioritizes older packets for transmission, is~optimal.

\section{Asymptotic Optimality of Index~Policies}
\label{s:aowip}
Although index policies are generally suboptimal for an RMABP where $M$ projects out of $N$ are selected at each time, a substantial amount of research has sought to establish the asymptotic optimality of index policies, most notably Whittle's, in a regime where both $M$ and $N$ increase to infinity in a fixed ratio. Assuming a homogeneous population of projects, Weber and Weiss~\cite{wewe90} demonstrated that Whittle's relaxation of the RMABP is asymptotically tight, in that the optimal average reward per project is the same for the relaxed and original problems in the aforementioned limit. They further showed that, although Whittle's conjecture on the asymptotic optimality of his index policy does not hold generally,  it does hold under a certain condition: when ``the differential equation describing the fluid approximation to the index policy has a globally stable equilibrium point.'' Although that condition is typically hard to check, they demonstrated in~\cite{wewe91} that it is satisfied in the case of three-state projects.

Bagheri and Scaglione~\cite{bagScag15} introduced a significant extension to the dynamic multi-channel access RMABP model with two-state channels, the so-called cognitive compressive sensing problem,  where the maximum number of channels to be sensed at each time can vary. They provided  conditions under which the myopic index policy is asymptotically~optimal. 

Larra\~{n}aga~et~al.~\cite{larraetal15} considered an RMABP model for the optimal scheduling of a multi-class queue with convex holding costs and user impatience under the average cost criterion. They derived index policies and demonstrated that Whittle's index policy is asymptotically optimal in both light- and heavy-traffic regimes. 

Ouyang~et~al.~\cite{ouyangetal16} considered a downlink scheduling problem modeled as a  POMDP RMABP problem. They established the asymptotic optimality of Whittle's index policy for two classes of positively correlated channels under the two-state channel model, assuming a recurrence condition that can be verified numerically. 

Verloop~\cite{verloop16} considered a set of priority policies for RMABPs with possibly non-indexable projects that can incorporate project arrivals and departures, as well as multi-action projects, which are asymptotically optimal when the differential equation for the system's fluid approximation
has a global attractor, similarly as in~\cite{wewe90}. She further demonstrated that these results can be applied to Whittle’s index policy in the case of indexable projects. 

Fu~et~al.~\cite{fuetal16} considered an RMABP dynamic job assignment model for a server farm with multiple heterogenous servers to optimize  energy efficiency. They showed that under certain conditions, Whittle’s index policy is asymptotically optimal, which requires a significant extension to the approach in~\cite{wewe90}. Motivated by geographically distributed server farms where available servers for a given job are job-dependent, Fu and Moran~\cite{fuMoran20} extended the model in~\cite{fuetal16}. They substantially improved on previous work by establishing not only the asymptotic optimality of an index policy but also its exponential convergence rate.  

Hu and Frazier~\cite{huFrazier17} considered a finite-horizon RMABP, derived an index policy based on optimal solutions of single-project subproblems, and argued its asymptotic optimality, showing an optimality gap of $o(N)$ for an $N$-project~model. 

Zayas-Cab\'an~et~al.~\cite{zayascetal19} considered a finite-horizon RMABP with multi-action projects and time-dependent upper bounds on the number of projects that can be selected at each time. They derived a heuristic policy based on the optimal solution to an LP relaxation and proved its asymptotic optimality. Their analysis does not rely on indexability or  stability conditions and applies to the model variant with project arrivals. The proposed policy was shown to have an optimality gap of $O(\sqrt{N} \log N)$.

Maatouk~et~al.~\cite{maatouketal21} considered an RMABP model for the optimal scheduling of transmissions over unreliable channels to minimize the average age of an information metric and establish the asymptotic optimality of the Whittle index policy under a recurrence condition that can be verified numerically. Kriouile~et~al.~\cite{kriouileetal22} considered a model extension and presented a novel approach to show the asymptotic optimality under less stringent~conditions. 

Brown and Smith~\cite{brownSmith20} presented index policies based on Lagrangian relaxation with an $O(\sqrt{N})$ optimality~gap.

In a promising recent work, Zhang and Frazier~\cite{zhangFrazier21} considered the finite-horizon RMABP and identified a non-degeneracy condition and a class of so-called  fluid-priority policies for which the asymptotic optimality gap was $O(1)$. When the condition failed to hold, they showed that fluid-priority policies still had an optimality gap of $O(\sqrt{N})$.

In a more recent work, Gast~et~al.~\cite{gastetal22} presented a framework for analyzing policies for both finite- and infinite-horizon RMABPs. Most notably, they provided conditions for a policy to be asymptotically optimal with an exponential convergence rate and presented a so-called LP-index policy, with~provably strong asymptotic optimality~properties. 

\section{Multi-Action~Bandits}
\label{s:marb}
In the standard RMABP introduced by Whittle~\cite{whit88b}, individual projects only allow two modes of operation, active and passive. However, a more general model was proposed much earlier, where individual projects were modeled as multi-action MDPs, which contained a passive action.
The resulting model, in which only one project can be active at each time (i.e., operated with an action different from the passive one) is known in the literature as a \emph{bandit superprocess} and~its inception was credited to Nash~\cite{nash73} by Gittins~\cite{gi79}.

Whittle~\cite{whittle80} provided an optimality proof for the Gittins index policy that was extended to bandit superprocesses, subject to a condition on dominating policies. The proof used a key construct later termed the Whittle integral. Brown and Smith~\cite{brownSmith13} showed that this integral
gives an upper bound on the value of a bandit superprocess, which is tighter than that obtained through Lagrangian relaxation.  They further showed how to efficiently compute the integral.
Hadfield-Menell and Russell~\cite{hadfieldetal15} also considered bandit superprocesses, providing a constructive definition of the Whittle integral and providing an alternate computation~method. 

The extension of the Whittle index to multi-action projects was first outlined by Weber \cite{weber07}, who illustrated it in a particular model and further outlined a means of computing the resulting index. 
In~\cite{nmvaluet08a}, the author formalized  the index extension in a general setting and demonstrated its applicability to a multi-armed multi-mode restless bandit problem, concerning the optimal dynamic allocation of a shared resource to a set of projects that can be operated in multiple modes while respecting a peak resource consumption limit. Sufficient PCL-based conditions were proposed to ensure both the existence of the index and the validity of an adaptive-greedy algorithm for its computation. 
These indices and further examples of their applications were also considered in the posterior work of Glazebrook~et~al.~\cite{glazetal11}.

The approach in~\cite{nmvaluet08a} was extended in~\cite{nmnetgcoop12} to a real-state project setting and demonstrated in a model for optimal dynamic energy management in a wireless sensor network. 
The approach outlined in~\cite{nmvaluet08a} was developed and established rigorously by the author in~\cite{nmmath22a}, which extended the sufficient indexability conditions in~\cite{nmaap01} for binary-action restless~bandits to \emph{multi-gear bandits}.

In recent years there has been an increased interest in RMABPs for multi-action projects, as they allow modeling of more realistic situations with projects that allow multiple operating modes. However, the heuristic policies that have been proposed are based on different approaches than the aforementioned indices. Zayas-Cab\'an~et~al.~\cite{zayascetal19} considered a policy for multi-action RMABPs obtained from an LP relaxation for which they established asymptotic optimality. The priority policies proposed by Verloop~\cite{verloop16} can also be applied to multi-action RMABPs. More recently, Killian~et~al.~\cite{killianetal21a,killianetal21} argued the relevance and importance of investigating multi-action RMABPs and presented powerful new methods for obtaining asymptotically optimal policies.
Xiong~et~al.~\cite{xiongetal22b,xiongetal22c} also considered multi-action RMABPs in a learning setting and developed efficient heuristic~policies.

\section{Lagrangian Index and Fluid Relaxation~Policies}
\label{s:lagpol}
Recent works have extended the Lagrangian relaxation approach used by Whittle~\cite{whit88b} to consider tighter Lagrangian relaxations with several multipliers, one per time period (for finite-horizon problems), as well as relaxations based on different ideas, in particular, fluid relaxations. These relaxations yield new policies with promising performance gains. 

Brown and Smith~\cite{brownSmith20} considered a  finite-horizon RMABP model with time-dependent rewards for dynamic item selection, including, for example, the problem of dynamic product assortment with demand learning by a retailer considered by Caro and Gallien~\cite{carogall07}. They proposed index policies and bounds based on a Lagrangian relaxation of the standard formulation that uses one Lagrange multiplier per period, and, under certain conditions, established their asymptotic optimality.

Brown and Zhang~\cite{brownZhang22a} studied an RMABP model, where  projects can consume different amounts of a shared resource, and there is exogenous information, modeled as a finite-state Markov chain, which can affect the shared resource limit, as well as each project’s rewards, transitions, and~resource consumption. They developed a Lagrangian relaxation and a ``dynamic fluid relaxation'' that provided upper bounds on the optimal value, as well as heuristic policies. The dynamic fluid relaxation bound and policy were shown to be asymptotically optimal, whereas those obtained from the Lagrangian relaxation were shown not to be asymptotically optimal.

Hao~et~al.~\cite{haoetal22} investigated a deadline scheduling problem with randomly arriving jobs that need to be processed before their deadlines expire, motivated by the problem of scheduling multiple electric vehicles in a charging station. They proposed a Lagrangian-based index policy, established its asymptotic optimality, and demonstrated through numerical experiments that the policy substantially outperforms  Whittle's index policy.

\section{Reinforcement Learning and Q-Learning~Approaches}
\label{s:reinflearn}
A difficulty that arises when addressing a real-world problem via a Markovian RMABP model is that its parameters (rewards and transition probabilities) are unknown, which has motivated research on this issue, among~which reinforcement learning and, in~particular, Q-learning (see Watkins and Dayan~\cite{watkinsDayan92}) approaches stand out (see, e.g.,~Chapters 6 and 7 in Bertsekas~\cite{bertsek12} and~the monograph by Powell~\cite{powell22}).

Fu~et~al.~\cite{fuetal19} considered the infinite-horizon average-cost RMABP and developed a tractable reinforcement learning algorithm based on parallel Q-learning recursions, which learns an approximation to the Whittle index. The approach was tested numerically and it was shown that its performance was close to that of Whittle's index policy in the model with all parameters known. 

Wu~et~al.~\cite{wuetal21} addressed the RMABP using a deep learning state-aware value function approximation approach, where the joint value function was approximated by a linear combination of individual project value~functions. 

Li~et~al.~\cite{lietal21} investigated a particular application of the RMABP to resource scheduling  in autonomous driving. They considered a Whittle index policy solution approach and approximated the index through a deep reinforcement learning method. 

Biswas~et~al.~\cite{biswasetal21a} considered RMABP models in public health settings and proposed Whittle-index-based Q-learning schemes that were  shown to converge to the performance of the Whittle index policy with all parameters~known.

Avrachenkov and Borkar~\cite{avraBorkar22} presented a reinforcement learning method  for RMABPs under the average-reward criterion and leveraged the properties of Whittle's index policy to reduce the search space of Q-learning and improve computational efficiency. They provided a convergence analysis and reported on numerical experiments that supported their findings.

Killian~et~al.~\cite{killianetal21} developed learning algorithms for multi-action RMABPs and combined  Lagrangian relaxation and Q-learning. They showed that under certain conditions, the proposed learning scheme converged to the asymptotically optimal multi-action RMAB policy. They further proposed another scheme that attains the asymptotic optimality properties of a Lagrange policy for multi-action RMABs through Q-learning. 

Nakhleh~et~al.~\cite{nakhlehetal21}  proposed a neural Whittle index network that can learn the Whittle indices for a given model by exploiting their properties. This motivated the use of deep reinforcement learning for training the neural network. The approach was demonstrated in computational experiments.

Nakhleh and Hou~\cite{nakhlehHou22} considered problems with optimal threshold policies to learn the optimal threshold. They developed an online policy that was shown to outperform other reinforcement learning algorithms by exploiting special structure. They further applied the results to efficiently learn the Whittle index.

\section{Regret-Based Online~Learning}
\label{s:learning}
In a significant area of research that considers finite-horizon bandit models with unknown parameters from an online learning perspective, the performance of the proposed policies is assessed using a measure of \emph{regret}, which compares the performance of a policy against that of an oracle, the~offline optimum, which implements the optimal policy based on known parameter values. The pioneering work in this vein was the paper by Robbins~\cite{robbins52}, who showed the existence of an asymptotically optimal policy, namely having a regret that grows sublinearly in the horizon $T$ for a two-armed Bernoulli bandit model with unknown parameters. In the MABP setting, Lai and Robbins~\cite{laiRobbins85} proved a landmark result, establishing a logarithmic  lower bound on the regret  for projects with i.i.d.\ rewards, as well as a policy asymptotically attaining it. Anantharam~et~al.~\cite{anantetal87} extended the results to MABPs with project transition probabilities parameterized by a~scalar.

More recent work has sought to extend these results to online learning in RMABP models with unknown parameters (see, e.g.,~Chapters 4 and 5 in the monograph by Zhao~\cite{zhao20}). As pointed out by Zhao, a key issue hindering regret analysis in such a setting is that the offline optimum when all parameters are known is typically unavailable, which is typically handled by considering instead  a proxy oracle that corresponds to a particular policy, leading to the notion of \emph{weak regret}. 

Major research efforts have been devoted to devising policies for these RMABPs, which are often not   index-based, along with corresponding regret analyses. It has been shown that, under suitable assumptions, some policies  attain a logarithmic regret for particular RMABP models, as demonstrated in the work of  Filippi~et~al.~\cite{filippietal11} and  Tekin and Liu~\cite{tekinLiu12}. A~more general regret bound  of order $O(\sqrt{T})$ was established by Ortner~et~al.~\cite{ortneretal14}. A different but related problem was studied by Garivier and Moulines~\cite{garivieretMoul11}, where project reward distributions were shown to abruptly change at unknown times.

Gupta~et~al.~\cite{guptaetal11} considered online RMABPs where unknown reward probabilities are assumed to drift over time and proposed a dynamic extension of Thompson sampling, which was shown to outperform alternative~methods. 

 Dai~et~al.~\cite{daietal11} addressed the online RMABP using an approach based on the assumption that the optimal policy belongs to a finite set of policies. They proposed a corresponding learning policy and illustrated it in an opportunistic spectrum access model. They established that the proposed  policy attains near-logarithmic regret. Liu~et~al.~\cite{liuetal13} proposed a different policy with logarithmic-order regret.  
 
Modi~et~al.~\cite{modietal17} focused on a particular online model of opportunistic spectrum access for which they provided an online learning policy with logarithmic-order regret.
    
Gr\"unew\"alder and Khaleghi~\cite{gruneKhal19} considered a more general setting where project dynamics have long-range dependence, including Markov chain dynamics as a special case. They proposed policies for this setting and provided corresponding regret analyses.
 
A policy with logarithmic-order regret was also obtained by Agrawal and Asawa~\cite{aggraAsawa20} in an opportunistic spectrum access model. Gafni and Cohen~\cite{gafniCohen21} considered a different class of RMAB problems motivated by communication networks and financial investment applications. They proposed a policy for which they established a logarithmic-order regret, as well as a finite-sample regret bound.

 Xu~et~al.~\cite{xuetal21} considered an online risk-averse RMABP model. They  proposed an index policy for the problem and showed that it attains a regret of order $O(\log T / T)$.

Gafni~et~al.~\cite{gafnietal22} provided an extension to the online RMABP that incorporated an exogenous global Markov process governing the rewards distribution of each project. They proposed a policy that achieved a logarithmic-order regret over time. 

Gafni and Cohen~\cite{gafniCohen22} considered a dynamic multi-channel access model for wireless networks, in which each channel has a different rate for each user. They proposed a policy for this model that attains a logarithmic-order regret. 
 
Xiong~et~al.~\cite{xiongetal22b} considered a finite-horizon online multi-action RMABP model. They proposed an index policy that is defined for non-indexable projects and showed that it achieves a sublinear regret with low computational complexity. The results were extended by 
Xiong~et~al.~\cite{xiongetal22c} to the average-reward~criterion.

\section{Applications: MDP~Models}
\label{s:appl}
In Ref.~\cite{whit88b}, Whittle illustrated the RMABP by sketching several potential applications, including the allocation of alternative medical treatments to patients for a condition caused by a continuously mutating virus, tracking surveillance of enemy submarines by aircraft, and dynamic activation and deactivation of a pool of workers subject to tiring and recovery. In Chapter 14.6 of his monograph~\cite{whit96}, he further outlined a machine maintenance model and derived its Whittle index.  
Numerous researchers have since demonstrated that a wide range of problems encountered in real-world applications can be modeled as RMABPs or variants thereof. This allows for the use of general solution methods developed for this framework. Below is a selection of such applications, focusing on RMABP models that can be formulated as MDPs with known model parameters.

\subsection{Variants of the~MABP}
\label{s:vmabp}
Several variants of the classic MABP can be reformulated as RMABPs. Consider, e.g.,~the MABP with switching penalties (costs or delays) considered by Banks and Sundaram~\cite{banksun94}, who showed that index policies are generally suboptimal in such a case, and~Asawa and Teneketzis~\cite{asatene96}, who proposed a heuristic index policy. In~\cite{nmvaluet207,nmijoc08,nmmath21}, the author demonstrated that the Asawa and Teneketzis index is the Whittle index of the problem in its restless reformulation and developed Whittle index algorithms for its efficient computation exploiting special structure. The results were extended to restless bandits with switching penalties, as outlined in~\cite{nmdags05}. Different policies for bandits with switching penalties were considered by Le Ny and Feron~\cite{lenyferon06}, Caro and Gallien~\cite{carogall07}, and~Arlotto~et~al.~\cite{arlottoetal14}. 

In~\cite{nmcdc05,nmijoc11}, the author considered the finite-horizon MABP and its extensions and showed that such problems can be reformulated as infinite-horizon RMABP models. This yielded a new efficient algorithm for computing a classic finite-horizon priority index by deploying a Whittle index algorithm in such a setting. Alternative finite-horizon index policies for the MABP were considered by Caro and Gallien~\cite{carogall07}. 

Dayanik~et~al.~\cite{dayaniketal08} considered a MABP where projects are not always available for selection. They showed that index policies are not optimal for this problem variant. They reformulated it as an RMABP and derived and analyzed its Whittle index policy. 

Caro and Yoo~\cite{caroYoo10} proposed a MABP model with random response delays. They showed that the resulting RMABP reformulation is indexable and computed the Whittle index for the special Beta-Bernoulli Bayesian learning model. Computational experiments showed that the index policy achieved near-optimal performance.

\subsection{Queueing~Models}
\label{s:qm}
Veatch and Wein~\cite{vewe96} formulated the problem of optimal scheduling in a multi-class M/M/1 make-to-stock queue as an RMABP model. They showed that the lost sales case is indexable and evaluated and tested the resulting Whittle index policy.  They further argued that the backorder case is non-indexable, which was also pointed out by Whittle in (\cite{whit96} Chapter~14.7). In \cite{nmmor06}, the author considered a model extension to a multi-class make-to-order/make-to-stock M/G/1 queue with possibly nonlinear holding costs in the backorder case and~showed that, in contrast to the aforementioned non-indexability results, the~resulting model was, in fact, indexable under an extended notion of indexability, relative to a mixed average-bias criterion (see also Ansell~et~al.~\cite{anselletal03} for an application of Whittle's index policy to the optimal scheduling of a make-to-order multi-class M/M/1 queue with nonlinear~costs).

In \cite{nmmp02}, the author illustrated the applicaiton of the PCL-based theoretical and algorithmic framework presented in that paper in the context of a general model for the optimal control of admission to a birth--death queue, which was solved by an extended Whittle index policy. The model was shown to be a building block for an RMABP formulation of a comprehensive model for the optimal dynamic control of admission  and routing to parallel queues. This formulation included features such as job abandonments and finite buffers and was further developed in~\cite{nmnetcoop07,nmcor12,nmejor12,nmcor19}.

Raissi-Dehkordi and Baras~\cite{raissi02} considered a model for optimal pull broadcast scheduling in information delivery systems, which they modeled through discrete-time bulk service queues. They derived and tested the Whittle index policy for an RMABP formulation of the model.

Dusonchet and Hongler~\cite{duhon03a}  derived the Whittle index for a make-to-stock queue with backorders under the discounted cost criterion. They further outlined a theory of Whittle indexability for continuous time and state restless bandits and demonstrated their approach using several random dynamic models, including diffusion processes.

Goyal~et~al.~\cite{goyaletal06} considered a discrete-time queueing model for the optimal scheduling of multimedia transmissions over a polled multiaccess fading channel, and~investigated the Whittle index policy for the resulting RMABP~formulation. 

In Ref.\ \cite{nmqs06}, the author formulated the problem of the optimal scheduling of a multi-class queue with finite buffers as an RMABP, established its indexability relative to a bias optimality criterion, and developed the corresponding Whittle index policy, which yielded nontrivial~insights.

Cao and Nyberg~\cite{caonyb08} proposed a Markovian model for the optimal dynamic admission control of multi-class traffic to a finite shared buffer. They established the model's indexability, evaluated the Whittle index, and numerically tested the resulting index policy. 

\textls[-20]{Borkar and Pattathil~\cite{borkarPatta22} investigated the so-called egalitarian processor sharing queueing model, which they reformulated as an RMABP and established its indexability. They showed how to compute the index and demonstrated its near-optimal performance through numerical experiments. }

\subsection{Web~Crawling}
\label{s:wc}
O'Meara and Patel~\cite{omearaPatel01} considered the problem of optimal scheduling of a web robot to  construct and maintain topic-specific web indexes. They realized that this problem could be modeled as an RMABP and developed a reinforcement learning algorithm for its approximate solution via an index policy similar to~Whittle's. 

In Ref.\ ~\cite{nmasmta14}, the author proposed a Markovian model for the optimal dynamic scheduling of page refreshes in a local repository of copies of randomly changing remote web pages. The model was reformulated as an RMABP and~Whittle's index policy was derived and tested. Avrachenkov and Borkar~\cite{avraBorkar18} considered a model for the optimal scheduling of a web crawler to retrieve ephemeral content from various sites. They showed that this model fits into the RMABP framework and developed a Whittle index policy, which they tested numerically.

\subsection{Public Health~Interventions}
\label{s:publichealth}
Deo~et~al.~\cite{deoetal13} developed a  model for optimal community-based healthcare delivery for a chronic disease, which was formulated as a variant of the finite-horizon RMABP.  They designed a myopic index heuristic policy and tested its performance using real data, demonstrating significant performance gains over the benchmark policy. 

Ayer~et~al.~\cite{ayeretal19} proposed an RMABP model to support the prioritization of hepatitis C treatment  decisions in U.S. prisons. They established indexability and  derived the Whittle index in closed form, deriving insights from it. They further proposed an adjusted closed-form index policy that was designed to overcome the limitations of Whittle's index policy. The model was validated using real-world data against benchmark policies.

Mate~et~al.~\cite{mateetal21,mateetal22} discussed the use of the RMABP framework in public health and reported on a field study aimed at assisting local health delivery agents to improve maternal and child health (see also the related work of Biswas~et~al.~\cite{biswasetal21a}).

\subsection{Communication~Networks}
\label{s:commnets}
Restless bandits have been widely deployed as a modeling framework in communication networks, which is currently one of its major application areas. In this type of setting, Wei~et~al.~\cite{weietal10} considered the problem of optimal relay selection in wireless cooperative networks, incorporating finite-state Markov channels, adaptive modulation and coding, and~residual relay energy. They showed that this problem can be formulated as an indexable RMABP.  Simulation results demonstrated  the effectiveness of the Whittle index~rule.

Wei and Neely~\cite{weiNeely16} proposed a model of power-aware throughput maximization in a multi-user file-downloading system. The model was formulated as a variant of the RMABP. An index policy that was different from Whittle's and based on a Lyapunov indexing approach was proposed and tested.

Aalto~et~al.~\cite{aaltoetal16} investigated optimal opportunistic scheduling  for downlink data traffic in a wireless cell with time-varying channels to minimize flow-level holding costs. They developed a size-aware index policy by deploying the Whittle index in a novel way. A numerical study demonstrated the improved performance achieved by the proposed policy.  

Borkar~et~al.\  \cite{borkaretal18} proposed a multi-user energy-efficient scheduling model in which each user has a separate queue and a cost is incurred for holding packets in each queue. Packets are transmitted through a shared channel with time-varying quality that can differ across users. Additionally, the cost incurred, i.e., energy consumed, for packet transmissions is a function of the channel quality. Indexability was proven for the average-cost criterion, Whittle's index was evaluated, and the resulting policy was tested. The Whittle index policy was shown to outperform previously considered policies such as max-weight scheduling and weighted fair scheduling.

Sun~et~al.~\cite{sunetal18} considered a model for heterogeneous cellular networks  with a macro base station and multiple small base stations. User equipment cell association was investigated to maximize the long-run average system throughput. The model was formulated as an RMABP for which index policies were derived and tested.

Aalto~et~al.~\cite{aaltoetal19} investigated the  opportunistic scheduling of downlink data traffic with partial channel information to minimize flow-level holding costs. The paper extended earlier work and, in particular, showed the indexability of  the  flow-level opportunistic scheduling problem with partial channel information in part of the parameter space. The authors derived a formula for the Whittle index, established optimality of threshold policies, and  tested numerically Whittle's index policy by comparing it against alternative policies.  

Wang~et~al.~\cite{wangetal19} considered a scheduling problem where a server opportunistically serves multiple user classes under time-varying multiple-state Markov channels with the goal of minimizing the average waiting cost. They reformulated the problem as an RMABP but noted that indexability was still open in this setting. They presented sufficient conditions on a channel-state transition matrix implying indexability and obtained the Whittle index in closed form. For the general case, they proposed an approximate Whittle~index. 

Sun~et~al.~\cite{sunetal20} addressed the Age of Information (AoI) minimization scheduling problem over a tar-topology wireless network, which was formulated as an RMABP. Indexability was shown and the Whittle index was obtained in closed form. The index was extended to incorporate more realistic features. Numerical studies demonstrated the effectiveness of the proposed index policies.

Chen~et~al.~\cite{chenetal22} proposed to use the uncertainty of information, as measured by Shannon’s entropy, as an information freshness metric. They considered a model where a central monitor observes $N$ binary Markov chains through $M<N$ communication channels and developed scheduling policies to minimize the long-run average uncertainty of information. The problem was cast as an RMAB and an index policy was developed and tested achieving excellent empirical performance. 

Singh~et~al.~\cite{singhetal22} considered the problem of determining which base station a user should associate with in a dense millimeter-wave network. The objective was to design an association policy to minimize the weighted average sojourn time for users in the system. The problem was formulated as an RMABP, which was shown to be  indexable. Whittle's index policy was derived and tested in a simulation study, in which it outperformed alternative user association policies.

\subsection{Miscellaneous~Applications}
\label{s:misc}
Given the huge variety of applications in which the RMAB framework has been deployed, this section discusses several papers that do not neatly fit into the other categories discussed previously. 

Caro and Gallien~\cite{carogall07} considered a model for the optimal dynamic product assortment of a retailer to maximize the overall profit for a selling season. The problem was formulated as a finite-horizon MABP with Bayesian learning in which multiple projects can be selected at each time. A closed-form index policy, which approximates Whittle's index policy, was derived, analyzed, and tested. 

Huberman and Wu~\cite{hubWu08} considered a model for the automatic generation of a ranking of 
information sources to be presented to users with limited attention. The objective was to maximize the total expected utility. The problem was formulated as an RMABP with dual-speed projects, which were shown to be PCL-indexable by Glazebrook~et~al.~\cite{gnma02}, and hence their Whittle index can be computed using the adaptive-greedy algorithm in~\cite{nmaap01}. 

Kumar and Saranga~\cite{kuSa10} developed near-optimal obsolescence mitigation policies based on an RMABP model and Whittle's index policy.

Temple and Frazzoli~\cite{tempFraz10} considered a well-studied online search problem known as the Cow
Path Problem (CPP), which is typically addressed via competitive analysis. They instead considered an MDP formulation for which they proved that the relaxed RMABP version was indexable, which allowed them to derive a Whittle index policy.

He~et~al.~\cite{heetal12} considered a model for the opportunistic scheduling of low-priority jobs onto under-utilized cloud resources left by high-priority jobs. Assuming  that the availability of servers to low-priority jobs can be modeled as on/off Markov chains, they formulated the problem as an RMAB. They established indexability and derived closed-form Whittle index formulae. A numerical study using real data center traces demonstrated the effectiveness of the policy compared with alternative rules.

Taylor and Mathieu~\cite{taylorMat14} proposed a dynamic demand response model formulated as an RMABP. Whittle  index policies were derived and used to rank loads. Numerical experiments showed that the resulting policy significantly outperforms the na\"{\i}ve greedy~policy. 

Sun and Ma~\cite{sunMa14} considered a crowd-sensing model in mobile social networks (MSNs) to maximize social welfare under a coverage constraint. The model was formulated as an RMABP, its indexability was proven, and~its Whittle index was used to design novel incentive schemes, which were shown to outperform previously proposed policies.

Lin~et~al.~\cite{linetal15} proposed a model for optimal forward-looking experiential learning problems, which was formulated as an RMABP. Under certain assumptions, indexability was established and the Whittle index and other index policies were derived and tested through numerical experiments with real data, which demonstrated its near-optimal utility and showed that they outperformed alternative~rules. 

Guo~et~al.~\cite{guoetal18} developed scheduling policies for networked cyber-physical systems that satisfy inter-delivery time requirements of clients connected through wireless channels. The problem was formulated as an infinite-state risk-sensitive MDP model. Among the results, when channels were not relatively reliable, they presented a  Whittle-like index policy for the model. Simulation results demonstrated the effectiveness of the proposed index policy. 

Yu~et~al.~\cite{yuetal18} addressed the stochastic deadline scheduling problem, which they formulated as an RMABP model.  This was shown to be indexable and the Whittle index was derived in closed form. The Whittle index policy was shown to be asymptotically optimal.

Qian~et~al.~\cite{qianetal16} considered a model for the optimal scheduling of patrol policies motivated by security domains with frequent interactions between defenders and attackers such as wildlife protection, which they formulated as an RMAB. They  provided sufficient conditions for indexability and an algorithm to test it. 

Avrachenkov~et~al.~\cite{avraetal17} investigated a model for Generalized Additive Increase Multiplicative Decrease (G-AIMD) dynamics for resource allocation under a fairness-based utility function. The model was formulated as an RMABP, indexability was established in special cases, and a means of computing the Whittle index was presented. The index policy was numerically tested and it was shown to achieve near-optimal performance.

Borkar~et~al.~\cite{borkarRavikSab17} considered an MDP model for optimal resource allocation in cloud computing based on dynamic pricing. The model was cast as an RMABP, indexability was proven, and an iterative scheme for computing the Whittle index was provided.

Motivated by the challenge of incorporating user feedback to tailor system operation for improved individual user satisfaction, Menner and Zeilinger~\cite{mennerZeil18} proposed a model to optimize the collection and processing of user feedback to  maximize a user comfort measure. The model was formulated as an RMABP, which was shown to be PCL-indexable. The Whittle index for this model was computed using the adaptive-greedy algorithm in~\cite{nmaap01}. Furthermore, the authors considered a learning model where transition probabilities are unknown and developed and tested an approach that combined restless bandit indices with upper confidence bound algorithms. 

Motivated by recommendation systems, Jhunjhunwala~et~al.~\cite{Jhunjetal18} considered an RMABP where, in each time period, the reward for selecting a project depends on the time that has elapsed since the project was last selected. They characterized the optimal policy with respect to the long-run average-reward criterion. 

Abbou and Makis~\cite{abbuMakis19} considered a maintenance planning model, where available repairmen are dynamically allocated to a set of unreliable production facilities with machines that incur losses due to degradation. The aim was to find a scheduling policy for maintenance interventions  that minimizes production losses per period. The model was formulated as an RMABP and was shown to be indexable.  The Whittle index policy was shown to achieve near-optimal performance and to outperform previously proposed policies.

Gerum~et~al.~\cite{gerumetal19} developed an RMABP model for  optimal inspection and maintenance scheduling policies in railway systems. Indexability was demonstrated and  Whittle indices were obtained, providing a novel index policy for such a system, which was shown to be highly effective in a data-driven setting. 

Li~et~al.~\cite{lietal20} considered an RMABP model for scheduling the allocation of limited resources to a large number of jobs such as medical treatments with random lifetimes and service times after a mass-casualty event, where jobs were initially subject to triage and classified. Whittle indices were derived, and another solution approach was considered through a nonstandard Lagrangian relaxation. Numerical experiments demonstrated that the second approach achieved better performance than the first one.

Fu and Moran~\cite{fuMoran20} studied a job-assignment model in a large-scale server farm system with geographically distributed heterogeneous servers. The goal was to maximize the system's energy efficiency through  dynamic  load control  on the networked servers. A scalable job-assignment policy was presented. Drawing on the asymptotic optimality analysis of Weber and Weiss~\cite{wewe90}, it was shown that the proposed policy quickly (exponentially) approaches asymptotic optimality, which was verified through numerical experiments. Fu~et~al.~\cite{fuetal22} extended the model and analysis to a setting with varying requests and limited-capacity resources that are shared by~requests. 

Dahiya et al.~\cite{dahiyaetal22} considered the dynamic allocation of human operators in a system with semi-autonomous robots that are required to perform independent task sequences, but are subject to the risk of getting stuck or failing. A human operator can assist the robot. The model was formulated as an RMABP and its indexability was proven under certain conditions. A simulation study demonstrated the near optimality of Whittle's index policy.

Motivated by mobile intervention problems, Ou~et~al.~\cite{ouetal22} utilized RMABs with network effects that were not amenable to computational solution using standard methods. They proposed a new solution approach for the networked RMABs, provided sufficient conditions for the optimality of their approach, and demonstrated its strong empirical performance in real-world scenarios.

\section{Applications: POMDP~Models}
\label{s:applPO}
A growing area of research considers RMABPs with continuous-state projects that arise from POMDP models (see, e.g.,~the monograph by Krishnamurthy~\cite{krishn16}), particularly in multi-target tracking and sensor scheduling applications (see also the monograph of Wang and Chen~\cite{wangChen21}).

La Scala and Moran~\cite{lascala06} considered an RMABP model for multi-target tracking with project states following scalar Kalman filter dynamics in a POMDP framework. Le Ny~et~al.~\cite{lenyetal08} and~Liu and Zhao~\cite{liuZhao08} proposed equivalent POMDP RMABP models that arise in different applications, the scheduling of unmanned aerial vehicles in the former paper and  multi-channel opportunistic access in the latter, with~imperfectly observed Gilbert--Elliott (two-state Markov chain) channels. Both papers established  indexability of the individual projects, evaluated Whittle's index policy, and tested it numerically, reporting excellent performance. 
 
 Liu and Zhao~\cite{liuZhao10} considered a class of RMABPs that arise in dynamic multi-channel access applications. They proved the indexability and derived the Whittle index in closed form for the
discounted- and average-reward criteria. The Whittle index policy was shown to be optimal in the case of identical projects under certain conditions. A numerical study demonstrated the near-optimal performance of Whittle's index policy and the related Lagrangian performance bound.

 Le Ny~et~al.~\cite{lenyetal11} considered the problem of scheduling observations using a number of sensors on a higher number of targets whose states obey continuous-time Kalman filter dynamics. Whittle's approach was deployed in this RMABP, indexability was argued, and in the scalar case with identical sensors, the Whittle index was derived in closed~form. 
 
Liu~et~al.~\cite{liuWebZhao11}  investigated a class of RMABPs where the active action resets the system's evolution. Indexability was established and the Whittle index was derived in closed form. The results were applied to opportunistic spectrum access and supervisory control systems. These results were extended by Akbarzadehand and Mahajan~\cite{akbarMahaj22b}.

Gan and Chen~\cite{ganChen12} introduced a new sensing policy for dynamic multi-channel access with primary and secondary users based on an RMABP formulation. Whittle's index policy was derived and applied and various results on the related Lagrangian relaxation were obtained.

He~et~al.~\cite{heetal11} considered a POMDP model  for topology tracking in a dynamic network with limited monitoring resources where links are modeled as independent on-off Markov chains. The goal was to maximize the overall  tracking accuracy of link states. A version of the model based on link sampling was formulated as an RMABP and its indexability was proven under certain conditions, allowing the deployment of Whittle's index policy for link sampling. Numerical results attested to the strong performance of the proposed approaches.

Meshram~et~al.~\cite{meshrametal15} considered a restless bandit model for a recommendation system. They characterized the optimal discounted policy for the single-project case with two underlying states and argued that in a certain special case the optimal policy is of threshold type, which they view as a relevant step towards establishing indexability in future work.

Ouyang~et~al.~\cite{ouyangetal15} considered a model of opportunistic multi-user scheduling in downlink networks with Markovian outage channels, which was cast as an RMABP. They showed that the model is indexable and obtained the Whittle index policy in closed form. Numerical experiments demonstrated the policy's near-optimal performance.

Taboada~et~al.~\cite{taboadaetal17} addressed the problem of scheduling traffic flows in wireless downlink systems under limited channel-state feedback to minimize the mean flow delay. The problem was cast into the RMABP framework with POMDP projects, indexability was argued, and~Whittle's index was evaluated and tested in numerical~experiments.

Meshram~et~al.~\cite{meshrametal18} considered an RMABP where each project can be in one of two states, which are not observable and need to be inferred from the current belief state and a possible binary signal. Single-project subproblems were shown to admit  an approximate threshold-type optimal policy in certain cases in which they satisfy an approximate indexability property. In cases where the optimality of threshold policies can be established, indexability was argued and the Whittle index was calculated. 

Elmaghraby~et~al.~\cite{elmaetal18} considered the problem of dynamic channel allocation for femtocells sharing the use of a regular macrocell spectrum. The channel state is not observable and macrocell user feedback  is utilized. The problem was cast as an RMABP with POMDP projects  and an approximate Whittle index policy was derived. A numerical study demonstrated the effectiveness of  the proposed policy compared to the myopic rule.

Mehta~et~al.~\cite{mehtaetal18} proposed an RMABP with constrained availability of projects and Markovian state evolution, with the true states being hidden. The optimality of a threshold policy was argued, and based on this, indexability was shown. A~formula for the Whittle index was derived for the rested case and an index algorithm was provided for the restless~case.

Kaza~et~al.~\cite{kazaetal19} considered a class of RMABPs with hidden states that allow cumulative feedback.  They showed that individual project subproblems, which they called lazy restless bandits (LRBs), have optimal policies of threshold type. Indexability was argued and the Whittle index was derived in closed form for two sets of special cases. An index-computing algorithm was provided and an extensive numerical study was reported.

Yang and Luo~\cite{yangLuo20} considered a massive multiple-input multiple-output (MIMO)  system with a smaller number of available orthogonal pilot sequences than users. The resulting pilot allocation problem was modeled as a POMDP with Gauss--Markov fading channels. The problem was cast as an RMABP and  the Whittle index was approximately derived. Numerical results were reported that demonstrated a strong performance for the index policy.

Hsu~et~al.~\cite{hsuetal20} considered a wireless broadcast network in which a base station updates users about random information arrivals subject to a transmission capacity constraint. The problem was cast into the RMABP framework with POMDP projects. Structural results on optimal policies were identified, allowing the deployment of Whittle's indexability and index policy. An online version was further considered. The results were validated in a numerical study. 

Wang~et~al.~\cite{wangetal20} studied dynamic channel allocation for the estimation of remote states in multi-agent systems. Whittle's approach was deployed in an RMABP formulation of the problem, and its strong performance  was demonstrated through numerical experiments.

Chen and Ephremides~\cite{chenEphre21} investigated a discrete-time model in which a base station simultaneously updates multiple users. The goal was to design a scheduling policy that minimizes an Age of Incorrect Information (AoII) metric for imperfect channel-state information. Whittle's index policy was derived under a simple condition. To avoid dealing with indexability, an alternative index policy was presented. A numerical study tested the performance of the policies  considered.

Kang and Joo~\cite{kangJoo21}  considered a model for minimizing  information mismatch under limited communication capabilities for a system in which a base station collects time-varying state information from multiple sources and makes decisions based on the collected information. An~RMABP formulation with POMDP projects was shown to be indexable and its Whittle index was obtained in closed form in certain~scenarios.

Motivated by public health interventions, Li and Varakantham~\cite{liVarak22} incorporated fairness constraints into an RMABP model. A~modified Whittle index was derived and, for the case where  transition probabilities are unknown, a~model-free learning method  was provided. The results were validated in a numerical~study.

\textls[-15]{Tong~et~al.~\cite{tongetal22} considered age-of-information minimization in Internet-of-Things networks with correlated sources. The problem was formulated as a correlated RMAB. A~generalized Whittle index and a generalized partial Whittle index were derived for the identical and non-identical channel settings, respectively, and corresponding index policies were proposed. A numerical study demonstrated that the policies are nearly-optimal, as they approach given Lagrangian-based bounds, and outperform state-of-the-art alternative~policies. }

\section{Conclusions}
\label{s:concl}
The work reviewed herein demonstrates without a doubt that after more than three decades since their inception, restless bandits continue to be a vibrant research area full of interesting problems to address in  theoretical, algorithmic, and application realms. It appears that the main avenues for further research should include the following: (1) deepening the understanding of indexability and developing easier-to-apply sufficient indexability conditions, both for models with binary-action and multi-action projects; (2) developing a full understanding of methods for designing asymptotically optimal index policies as the size of the model scales; (3) devising methods suitable for very large-scale models arising from applications; (4) developing and refining methods for designing and analyzing effective policies for models with unknown parameters that need to be learned online; and (5) further extending the scope and modeling power of the RMABP framework.


\vspace{6pt} 




\section*{Funding}
This work was funded in part by the Spanish State Research Agency (\emph{\hl{Agencia}  Estatal de Investigaci\'on}, AEI) under grant PID2019-109196GB-I00/AEI/10.13039/501100011033 and by the \emph{Comunidad de Madrid} in the context of a multi-year agreement with Carlos III University of Madrid within the activity ``\emph{Excelencia para el Profesorado Universitario}'' in~the framework of the V Regional Plan of Scientific Research and Technological Innovation 2016--2020.
 
\section*{Abbreviations}
{The following abbreviations are used in this manuscript:\\

\noindent 
\begin{tabular}{@{}ll}
DP & Dynamic programming \\
LP & Linear programming \\
MDP & Markov decision process\\
POMDP & Partially observable Markov decision process \\
MABP & Multi-armed bandit problem\\
RMABP & Restless multi-armed bandit problem\\
\end{tabular}
}

\end{document}